\documentclass[12pt,reqno,a4paper]{amsart}
\usepackage{
    amsmath,  amsfonts, amssymb,  amsthm,   amscd,
    gensymb,  graphicx, comment,  etoolbox, url,
    booktabs, stackrel, mathtools,enumitem, mathdots,  microtype, lmodern,    mathrsfs, graphicx, tikz,  longtable,tabularx, float, tikz, pst-node, tikz-cd, multirow, tabularx, amscd,  bm, array, makecell, diagbox, booktabs,ragged2e, caption, subcaption }
\usepackage{makecell,slashbox}

\usepackage{xcolor}
\usepackage[utf8]{inputenc}
\usepackage{microtype, fullpage, wrapfig,textcomp,mathrsfs,csquotes,fbb}
\usepackage[colorlinks=true, linkcolor=blue, citecolor=blue, urlcolor=blue, breaklinks=true]{hyperref}
\usepackage[capitalise]{cleveref}
\usepackage{todonotes}
\usetikzlibrary{positioning}
\usetikzlibrary{shapes,arrows.meta,calc}
\usetikzlibrary{arrows}

\newtheorem{theorem}{Theorem}[section]

\newtheorem{definition}[theorem]{Definition}
\newtheorem{example}[theorem]{Example}

\newtheorem{proposition}[theorem]{Proposition}
\newtheorem{remark}[theorem]{Remark}
\newtheorem{discussion}[theorem]{Discussion}

\newcommand\R{\mathbb{R}}
\newcommand\Z{\mathbb{Z}}
\newcommand\C{\mathbb{C} }

\newcommand{\TC}{\mathrm{TC}}

\newcommand{\ct}{\mathrm{cat}}
\newcommand{\cl}{\mathrm{cup}}
\newcommand{\zl}{\mathrm{zcl}}

\newcolumntype{x}[1]{>{\centering\arraybackslash}p{#1}}

\begin{document}
\title{LS-category and topological complexity of Milnor manifolds and corresponding generalized projective product spaces}
\author[N. Daundkar ]{Navnath Daundkar}
\address{Department of Mathematics, Indian Institute of Technology Bombay, India}
\email{navnathd@iitb.ac.in}
\begin{abstract}
Milnor manifolds are a class of certain codimension-$1$  submanifolds of the product of projective spaces. In this paper, we study the LS-category and topological complexity of these manifolds. We determine the exact value of the LS-category and, in many cases, the topological complexity of these manifolds. We also obtain tight bounds on the topological complexity of these manifolds. 
It is known that Milnor manifolds admit $\Z_2$ and circle actions. We compute bounds on the equivariant LS-category and equivariant topological complexity of these manifolds. Finally, we describe the mod-$2$ cohomology rings of some generalized projective product spaces corresponding to Milnor manifolds and use this information to compute the bound on LS-category and topological complexity of these spaces.
\end{abstract}
\maketitle

\section{Introduction}   \label{sec:intro}
Farber  \cite{FarberTC} introduced the notion of topological complexity to study the robot motion planning problem from a topological perspective. 
For a topological space $X$, this numerical homotopy invariant is denoted as $\TC(X)$. 
Let $X$ be a path connected space and $PX$ be the space of all paths
in $X$ equipped with a compact open topology. 
Let $\gamma:[0,1]\to X$ be any path in $X$. 
Then there is a fibration $\pi : PX \to X\times X$ defined by $\pi(\gamma)=(\gamma(0),\gamma(1))$.
The \emph{topological complexity} of $X$ is the smallest $k$ such that $X\times X$ admit an open cover $U_1, \dots, U_k$ such that there exist continuous sections of $\pi$ on $U_{i}$ for $1\leq i \leq k$. In general, determining the exact value of  $\TC(X)$ is a hard problem. 
Since Farber first introduced this invariant in $2003$, numerous methods have been created to determine bounds on $\TC(X)$, and the precise value of this invariant has been calculated for many spaces. 
This invariant is closely related to the well known invariant, the \emph{Lusternik-Schnirelmann category} (\emph{LS-category}). For a path connected space $X$, the LS-category of $X$ is denoted as $\mathrm{cat}(X)$. The $\mathrm{cat}(X)$ is the smallest integer $r$ such that $X$ can be covered by $r$ open subsets $V_1, \dots, V_r$ with inclusion $V_i\xhookrightarrow{} X$ is null-homotopic for $1\leq i\leq r$. 

In this paper, we study these invariants for Milnor manifolds.
In \cite{Milnor}, Milnor introduced the class of submanifolds of the products of real and complex projective spaces to define  generators for the unoriented cobordism algebra. 
Let $r$ and $s$ be integers such that $0\leq s\leq r$. A real Milnor manifold, denoted by $\mathbb{R}H_{r,s}$, is the submanifold of the product $\mathbb{R}P^r \times \mathbb{R}P^s$ of real projective spaces defined as follows:
\[\mathbb{R}H_{r,s}:=\Big\{\big([x_0: \dots: x_r], [y_0:\dots:y_s]\big) \in \mathbb{R}P^r \times \mathbb{R}P^s~|~ x_0y_0 +\hspace{1mm}\cdots \hspace{1mm}+ x_sy_s=0\Big\}.\] 
The real Milnor manifold $\mathbb{R}H_{r,s}$ is an $(s+r-1)$-dimensional closed, smooth manifold.
Observe that, we have the following fiber bundle 
\begin{equation}\label{eq: fbRHrs}
 \mathbb{R}P^{r-1} \stackrel{i}{\hookrightarrow} \mathbb{R}H_{r,s} \stackrel{p}{\longrightarrow} \mathbb{R}P^{s},
\end{equation}
where $p$ is induced from the projection onto $\R P^s$. More explicitly, \[p(([x_0: \dots: x_r], [y_0:\dots:y_s]))=[y_0,\dots,y_s].\]
The complex Milnor manifold, denoted by $\C H_{r,s}$ can be defined similarly
\[\C H_{r,s}:=\Big\{\big([z_0:\dots:z_r],[w_0:\dots:w_s]\big)\in \C P^r \times \C P^s~|~ z_0\overline{w}_0+\cdots+z_s\overline{w}_s=0\Big\}.\] 
Note that $\C H_{r,s}$ is a $2(s+r-1)$-dimensional closed, smooth manifold.
As in the real case, $\C H_{r,s}$ is the total space of the fiber bundle 
\begin{equation}\label{eq: fbCHrs}
 \C P^{r-1} \stackrel{i}{\hookrightarrow} \C H_{r,s} \stackrel{p}{\longrightarrow} \C P^{s},  
\end{equation}
where $p$ is induced from the projection onto $\C P^s$.
Note that, $\mathbb{R}H_{r,0}= \mathbb{R}P^{r-1}$ and $\C H_{r,0}= \C P^{r-1}$.
Milnor \cite[Lemma 1]{Milnor} shown that the unoriented cobordism algebra of smooth manifolds is
generated by the cobordism classes of real projective spaces and real Milnor manifolds. Dey and Singh \cite{DeyMilnormfds} characterized milnor manifolds which admit free $\Z_2$ and $S^1$-actions with restriction on $r$ and computed the equivarinat cohomology rings in this case. 

The purpose of this paper is to compute the LS-category and topological complexity of Milnor manifolds. We also study the generalized projective product spaces corresponding to Milnor manifolds (see \Cref{sec: cat and tc gpps}).

The paper is organized as follows: We begin \Cref{sec: prel} by recalling recall some basic notions in the theory of LS-category and topological complexity. Then we state some results about the (equivariant) cohomology rings of Milnor manifolds.

In \Cref{sec: cat and tc}, we compute the exact value of the LS-category of Milnor manifolds and obtain the sharp lower and upper bounds on the topological complexity of Milnor manifolds (see \Cref{prop: cat} and \Cref{thm: TClbub}). As an application of these results, we get the exact value of the topological complexity in some specific cases (see \Cref{ex: exact tc}). Then we showed that
the value of the modulus $|\TC(\R P^r)+\TC( \R P^s)-1-\TC(\R H_{r,s})|$ is either $0$ or $1$ for infinitely many values of $r$ and $s$.

In \Cref{sec: eq cat and tc}, we compute bounds on the equivariant LS-category and equivariant topological complexity of Milnor manifolds .

The \Cref{sec: cat and tc gpps} is devoted to the study of generalized projective product spaces whose one of the components is the Milnor manifold. We describe the mod-$2$ cohomology ring of these spaces and compute the bounds on their LS-category and the topological complexity.

\section{LS-category, topological complexity and Milnor manifolds}\label{sec: prel}
This section is divided into two parts. We begin the first part with some definitions. Then we state some results related to LS-category and topological complexity.
In the second part, we present some results related to the Milnor manifolds that we use to compute lower bounds on topological complexity.

\subsection{LS-category and topological complexity}
The sectional category of a Hurewicz fibration was introduced by \v{S}varc in  \cite{Sva}. Colman and Grant \cite{colmangranteqtc} generalized the notion of sectional category and introduced the equivariant sectional category. Let $G$ be a group and $E, B$ be $G$-spaces such that $p\colon E \to B$ is a $G$-map. The symbol $\simeq_G$ denotes the $G$-homotopy equivalence.
\begin{definition}[{\cite[Definition 4.1]{colmangranteqtc}}]
The equivariant sectional category of a map $p$, denoted by ${\rm secat}_G(p)$, is the smallest positive integer $k$ such that there is a $G$-invariant open cover $\{U_1, \ldots, U_k\}$ of $B$ and a $G$-map $\sigma_i \colon U_i \to E$ satisfying $p\circ \sigma_i \simeq_G \iota_{U_i}$ for $i=1, \ldots, k$, where $\iota_{U_i} \colon U_i \hookrightarrow B$ is the inclusion.
\end{definition}
If no such $k$ exists, we say ${\rm secat}_G(p)=\infty$.
If $G$ is a trivial group, then ${\rm secat}_G(p)$ is called the sectional category of $p$. Moreover, if $E$ is contractible, then the ${\rm secat}_G(p)$ coincided with the LS-category of the base space,  $B$.  
The topological complexity of $B$ is a sectional category of a free path space fibration $\pi: PB\to B\times B$, defined by $\pi(\gamma)=(\gamma(0),\gamma(1))$.

For $G$-spaces, Marzantowicz \cite{Eqlscategory} introduced the notion of an equivariant LS-category. This is a homotopy invariant of a $G$-space $B$, denoted by $\ct_G(B)$.
\begin{definition}
The equivariant LS-category $\mathrm{cat}_G(B)$ is the smallest integer $k$ such that $B$ can be covered by $k$ many $G$-invariant open subsets $U_1, \dots, U_k$ of $B$ with each inclusion $U_j\xhookrightarrow{} B$  $G$-homotopic to an orbit inclusion $Gb_j \hookrightarrow B$ for some $b_j \in B$ for $j=1, \ldots, k$.    
\end{definition}
The sets $U_j$'s are called $G$-categorical open subsets of $B$. We note that if $b \in B^G$ and $B$ is a $G$-connected space, then ${\rm secat}_G(\iota) = \ct_G(B)$ for the inclusion  map $\iota \colon \{b\} \to B$, see \cite[Corollary 4.7]{colmangranteqtc}. Also, if $G$ is trivial, then $\ct_G(B)=\ct(B)$ and $U_j$'s are called categorical open subsets of $B$. 

The following result by Farber shows how the LS-category and topological complexity are related. 
\begin{theorem}[{\cite[Theorem 4, Theorem 5]{FarberTC}}\label{thm: usual lbub}]
If  $X$ is a path connected, paracompact topological space. Then
\begin{equation}
 \ct(X)\leq \TC(X)\leq \ct(X\times X)\leq 2\dim(X)-1. 
\end{equation}   
\end{theorem}

Farber proved the product inequality of topological complexity.
\begin{theorem}[{\cite[Theorem 11]{FarberTC}}]
Let $X$ and $Y$ be path connected spaces. Then
 \begin{equation}\label{eq: tc_prod}
   \TC(X\times Y)\leq \TC(X)+\mathrm{TC}(Y)-1. 
\end{equation}   
\end{theorem}
Fox proved the product inequality for the LS-category.
\begin{theorem}[{\cite{Fox}}]
If $X$ and $Y$ are the path connected spaces and $X \times Y$ is completely normal
 \begin{equation}\label{eq_ls_cat_prod}
   \mathrm{cat}(X\times Y)\leq \mathrm{cat}(X)+\mathrm{cat}(Y)-1. 
\end{equation}   
\end{theorem}
We now describe algebraic methods to obtain lower bounds on these homotopy invariants.
Let $R$ be a commutative ring with unity. Let $H^{\ast}(X;R)$ be the cohomology ring of a path connected topological space $X$. The cup-length of $X$ over $R$ is the maximal number $k$ such that there exists $x_i\in H^{\ast}(X;R)$ for $i=1, \ldots, k$ with $\prod_{i=1}^{k}x_i\neq 0$. We denote this number by $\cl_{R}(X)$. It was shown in \cite{CLOT} that the number gives a lower bound for the $\ct(X)$.
\begin{theorem}[{\cite[Proposition 1.5]{CLOT}}]
Let $X$ be a path connected topological space. Then 
\begin{equation}\label{eq: cup-length}
\cl_{R}(X)+1 \leq \ct(X).    
\end{equation}  
\end{theorem}
The dimensional upper bound on the LS-category is given as follows.
\begin{theorem}[{\cite[Proposition 1.7]{CLOT}}\label{thm: catdimub}]
 Let $X$ be a path connected, locally contractible paracompact space. Then \[\ct(X)\leq \dim(X)+1,\] where $\dim(X)$ is the covering dimension of $X$.  
\end{theorem}

Let $\cup \colon H^{\ast}(X;R)\otimes H^{\ast}(X;R) \longrightarrow H^{\ast}(X;R)$ be the cup product map. Then the zero-divisor-cup-length of $X$ with respect to the coefficient ring $R$ is defined as the maximal number $k$ such that there exist elements $u_i\in H^{\ast}(X;R)\otimes H^{\ast}(X;R)$ such that $u_i\in \ker(\triangle^{\ast})$, where $\triangle$ is the diagonal map $X \to X \times X$
and $u_1\cup \dots \cup u_k\neq 0$. We denote this number by $\zl_{R}(X)$. Farber showed that the number $\zl_{R}(X)$ gives a lower bound for the $\TC(X)$. 
\begin{theorem}[{\cite[Theorem 7]{FarberTC}}]
 Let $X$ be a path connected topological space. Then  
 \begin{equation}\label{eq: zcl}
 \zl_{R}(X) +1\leq \TC(X).    
 \end{equation}
\end{theorem}

\subsection{Milnor manifolds}

The mod-$2$ cohomology ring of $\R H_{r,s}$ and $\C H_{r,s}$ was computed in \cite{Buchstaber} and \cite{Mukherjee}).

\begin{theorem}[{\cite[Lemma 3]{Mukherjee}}]\label{thm: coring}
Let $0 \leq s \leq r$. Then
\[H^{\ast}(\R H_{r,s}; \Z_2) \cong \displaystyle\frac{\Z_2[a,b]}{\langle a^{s+1}, b^r+ab^{r-1}+ \dots+a^sb^{r-s} \rangle},\]
where both $a$ and $b$ are homogeneous elements of degree one.

\end{theorem}

\begin{theorem}[{\cite{Mukherjee}}\label{thm: coring CHrs}]
Let $0 \leq s \leq r$. Then
\[H^{\ast}(\C H_{r,s}; \Z_2) \cong \Z_2[c,d]/\langle c^{s+1},h^r+cd^{r-1}+\dots+c^sd^{r-s} \rangle,\]
where both $c$ and $d$ are homogeneous elements of degree two.   
\end{theorem}

We list here some results from \cite{DeyMilnormfds} related to free actions of $\Z_2$ and $S^1$ on Milnor manifolds that we are going to use to compute lower bounds on the (equivariant) LS-category and topological complexity of Milnor manifolds.

\begin{theorem}[{\cite[Corollary 5.6, Corollary 5.8]{DeyMilnormfds}}\label{thm: existence of invo}\label{thm: free invo milnor mfds}]
Let $1<s<r$ such that $r\not\equiv2({\mod 4})$.
 Then $F H_{r,s}$ admit a free involution if and only if both $r$ and $s$ are odd integers, where $F$ is either $\R$ or $\C$ .
\end{theorem}

\begin{theorem}[{\cite[Proposition 5.9]{DeyMilnormfds}}\label{thm: free S1 action}]
Let $1\leq s\leq r$. Then $S^1$ acts freely on $\R H_{r,s}$ if and only if both $r$ and $s$ are  odd integers.   
\end{theorem}

\begin{theorem}[{\cite[Proposition 5.1, Proposition 5.2]{DeyMilnormfds}}\label{thm: freeZ2 action}]
Let $1<s<r$ and $r\not\equiv 2({\rm mod}~4)$ such that $F H_{r,s}$ admit a free $\Z_2$-action. Then the induced action on $H^*(F H_{r,s})$ is trivial.
\end{theorem}

\begin{theorem}[{\cite[Theorem 6.1]{DeyMilnormfds}}\label{thm: eq coho ring real milnor mfds}]
Suppose $\Z_2$ acts freely on $\R H_{r,s}$ such that it induces a trivial action on cohomology. Then     
\[H^{\ast}(\R H_{r,s}/\Z_2; \Z_2) \cong \Z_2[x,y,z,w]/I,\] where $|x|=|y|=2$, $|z|=|w|=1$ and the ideal $I$ is generated by following relations: \begin{enumerate}
    \item $z^2$, 
    \item $w^2-\gamma_1 zw - \gamma_2 x - \gamma_3 y$,
    \item $ x^{\frac{s+1}{2}} + \alpha_0zwx^{\frac{s-1}{2}} + \alpha_1zwx^{\frac{s-3}{2}}y+\cdots+ \alpha_{\frac{s-1}{2}}zwy^{\frac{s-1}{2}}$,
    \item $(w+\beta_0z)y^{\frac{r-1}{2}}+(w+\beta_1z)xy^{\frac{r-3}{2}}+\cdots+ (w+\beta_{\frac{s-1}{2}}z)x^{\frac{s-1}{2}}y^{\frac{r-s}{2}} $,
\end{enumerate}
where $\alpha_i, \beta_i, \gamma_i \in \mathbb{Z}_2$.
\end{theorem}

\begin{theorem}[{\cite[Theorem 6.2]{DeyMilnormfds}}\label{thm: eq coho ring complex milnor mfds}]
Suppose $\Z_2$ acts freely on $\C H_{r,s}$ such that it induces a trivial action on cohomology. Then     
\[H^{\ast}(\C H_{r,s}/\Z_2; \Z_2) \cong \Z_2[x,y,z,w]/I,\] where $|x|=|y|=4$, $|z|=1$, $|w|=2$ and the ideal $I$ is generated by following relations: \begin{enumerate}
    \item $z^3$, 
    \item $w^2-\gamma_1 z^2w - \gamma_2 x - \gamma_3 y, x^{\frac{s+1}{2}} + \alpha_0z^2wx^{\frac{s-1}{2}} + \alpha_1z^2wx^{\frac{s-3}{2}}y+\cdots+ \alpha_{\frac{s-1}{2}}z^2wy^{\frac{s-1}{2}}$,
    \item  $(w+\beta_0z^2)y^{\frac{r-1}{2}}+(w+\beta_1z^2)xy^{\frac{r-3}{2}}+\cdots+ (w+\beta_{\frac{s-1}{2}}z^2)x^{\frac{s-1}{2}}y^{\frac{r-s}{2}}$
\end{enumerate}
where $\alpha_i, \beta_i, \gamma_i \in \mathbb{Z}_2$.
\end{theorem}

\begin{theorem}[{\cite[Theorem 6.3]{DeyMilnormfds}}\label{thm: S1 eq coho ring}]
Suppose $S^1$ acts freely on $\R H_{r,s}$. Then     
\[H^{\ast}(\R H_{r,s}/S^1; \Z_2) \cong \Z_2[x,y,w]/I,\] where $|x|=|y|=2$, $|w|=1$ and the ideal $I$ is generated by following relations: \begin{enumerate} 
    \item $x^{\frac{s+1}{2}}$,
    \item $wy^{\frac{r-1}{2}}+xwy^{\frac{r-3}{2}}+\cdots+ wx^{\frac{s-1}{2}}y^{\frac{r-s}{2}}$,
    \item $w^2-\alpha x - \beta y$ .
\end{enumerate}    
\end{theorem}

\section{LS-category and topological complexity of Milnor manifolds}\label{sec: cat and tc}
 In this section, we use the mod-$2$ cohomology ring information of Milnor manifolds to compute bounds on  the LS-category and topological complexity of the real Milnor manifolds and complex Milnor manifolds. We compute the exact value of these invariants in infinitely many cases. We note that all the cohomological calculations are done with respect to $\Z_2$ coefficients.
\begin{theorem}\label{prop: cat}
Let $0 \leq s \leq r$. Then
 $\ct(F H_{r,s})=r+s$ where  $F$ be either $\R$ or $\C$. 
\end{theorem}
\begin{proof}
Since both  $\R H_{r,s}$ and $\C H_{r,s}$ are $\Z_2$-orientable, the product $a^sb^{r-1}$ and $f^sg^{r-1}$ of homogeneous classes of degree-$1$ and degree-$2$ are non-zero in $H^{s+r-1}(\R H_{r,s};\Z_2)$ and $H^{2(s+r-1)}(\C H_{r,s};\Z_2)$, respectively. Therefore, the  cup-length of both  $\R H_{r,s}$ and $\C H_{r,s}$ over $\Z_2$ is $s+r-1$. Then \eqref{eq: cup-length} 
give $s+r\leq \ct(\R H_{r,s})$ and $s+r\leq \ct(\C H_{r,s})$. Recall that $\R H_{r,s}$ is an $s+r-1$ dimensional smooth, closed manifold. Therefore, using \Cref{thm: catdimub}
we get $\ct(\R H_{r,s})\leq s+r$. We have a fibre bundle $\C P^{r-1} \stackrel{i}{\hookrightarrow} \C H_{r,s} \stackrel{p}{\longrightarrow} \C P^{s}$. This implies $\C H_{r,s}$ is simply connected, as both $\C P^{r-1}$ and $\C P^{s}$ are simply connected. Then \cite[Proposition 5.1]{James}, give $\ct(\C H_{r,s})\leq s+r$. This proves the proposition.
\end{proof}

\begin{theorem}\label{thm: TClbub}
Let $0 \leq s \leq r$ be two integers.
\begin{enumerate}
\item Suppose $F$ is either $\R$ or $\C$. Then
\begin{equation}\label{eq: TC Milmfds}
 s+r\leq \TC(F H_{r,s})\leq 2(s+r)-1.  
\end{equation} 
\item If either $r>2$ and $s>1$ or $r=1$ and $s=1$, then
\begin{equation}\label{eq: tcsharpub}
 \TC(\R H_{r,s})\leq 2(s+r-1).   
\end{equation}
\item Suppose $s=2^{t_1}+1$ and $r=2^{t_2}$ for some positive integers $t_1, t_2$. Then 
\begin{equation}\label{eq: equaTCMil}
 2(s+r-1)-1\leq \TC(F H_{r,s})\leq 2(s+r-1)+1,
\end{equation}
where $F$ is either $\R$ or $\C$.
Moreover, if $t_2>1$, then 
\begin{equation}\label{eq: equaTCMil1}
 2(s+r-1)-1\leq \TC(\R H_{r,s})\leq 2(s+r-1).
\end{equation}
\end{enumerate}
\end{theorem}
\begin{proof}
The inequalities of \eqref{eq: TC Milmfds} follow from \Cref{thm: usual lbub}.
\Cref{prop: cat}. 

Note that, if $r>2$ and $s>1$, then one can use the homotopy exact sequence of \eqref{eq: fbRHrs} to show that $\pi_1(\R H_{r,s})=\Z_2\times \Z_2$. 
Then part $(1)$ and part $(2)-2c$ of \cite[Theorem 1.2]{abfundamentaltc} implies that the $\TC(\R H_{r,s})$ is non-maximal. In particular, we get the inequality \eqref{eq: tcsharpub}.
If $(r,s)=(1,1)$, then it is clear that $\R H_{1,1}\cong S^1$. Thus \eqref{eq: tcsharpub} follows. 
This completes the proof of the inequality in \eqref{eq: tcsharpub}.

We now compute the zero-divisor-cup-length to show the left inequality in \eqref{eq: equaTCMil}. Let $\bar{a}=1\otimes a+a\otimes 1$ and $\bar{b}=1\otimes b+b\otimes 1$ be the two basic zero divisors. Then 
\begin{align*}
  \bar{a}^{2(s-1)-1}&=\sum_{i=0}^{2s-3}\binom{2s-3}{i} a^{2s-i-3}\otimes a^{i}  \\
  & =\sum_{j=0}^{3}\binom{2s-3}{s-j}a^{s-3+j}\otimes a^{s-j} \hspace{.5cm}(\text{since}~~ a^{s+j}=0, j\geq 1).
\end{align*}
Note that $b^{r+j}=0$ for $j\geq 1$. Otherwise, suppose $b^{r+1}\neq 0$. Since $a^sb^{r-1}$ generates $H^{r+s-1}(\R H_{r,s})$, $a^{s-2}b^{r+1}=a^{s}b^{r-1}$. From \Cref{thm: coring}, we have $b^r=ab^{r-1}+\dots+a^sb^{r-s}$. This implies $a^{s-2}(ab^{r-1}+\dots+a^sb^{r-s})b=a^sb^{r-1}$. Since $a^{s+j}=0$ for $j\geq 1$, we get $a^{s-1}b^r+a^sb^{r-1}=a^sb^{r-1} $. Consequently, we have $a^{s-1}b^r=0$. This is a contradiction to the fact that $a^{s-1}b^r$ generates the top cohomology of $\R H_{r,s}$. Now consider
\[\bar{b}^{2r-1}=\sum_{i=0}^{2r-1}\binom{2r-1}{i} b^{2r-i-1}\otimes b^{i}=\binom{2r-1}{r-1}b^r\otimes b^{r-1}+\binom{2r-1}{r}b^{r-1}\otimes b^r.\]
It is well known that $\binom{2^k-1}{l}$ is odd for all $l$. Therefore, for $s=2^{t_1}+1$ and $r=2^{t_2}$, both $  \bar{a}^{2(s-1)-1}$ and $\bar{b}^{2r-1}$ are non-zero. 
 Consider the following product 
 \[\bar{a}^{2(s-1)-1}\cdot \bar{b}^{2r-1}=\sum_{j=0}^{2}a^{s-3-j}b^r\otimes a^{s-j}b^{r-1}+\sum_{j=1}^{3}a^{s-3-j}b^{r-1}\otimes a^{s-j}b^{r}.\]
 It can be seen that a set $\{a^{s-j}b^r, a^{s-j+1}b^{r-1}\}$ is a part of the basis of $H^{r+s-j}(\R H_{r,s})$. Therefore, the product $\bar{a}^{2(s-1)-1}\cdot \bar{b}^{2r-1}$ is non-zero. Therefore, using \eqref{eq: zcl} we get \[2(s+r-1)-1\leq \TC(\R H_{r,s}).\]
Since the cohomology ring of $\C H_{r,s}$ has a similar description as that of $\R H_{r,s}$, we can do similar calculation to show that $2(s+r-1)-1\leq \TC(\C H_{r,s})$. 
 
 Note that the right inequality of \eqref{eq: equaTCMil} follows from the \Cref{prop: cat} and \Cref{thm: usual lbub}. We can use \Cref{thm: coring CHrs} and do similar calculations when $F=\C$. 
 The right inequality of \eqref{eq: equaTCMil1} follows form \eqref{eq: tcsharpub}.
\end{proof}

\begin{remark}
Note that, if either $s=2^{t_1}+1$ or $s=2^{t_1'}$ and $r=2^{t_2}$, then \[\bar{a}^{2s-1}\cdot \bar{b}^{2r-1}=a^{s-1}b^r\otimes a^sb^{r-1}+a^sb^{r-1}\otimes a^{s-1}b^r=0.\]
Therefore, the inequality in \eqref{eq: equaTCMil} is sharp.
\end{remark}

\begin{example}\label{ex: exact tc}
At this point, we make the following observations.
\begin{enumerate}
\item Suppose $(r,s)=(1,1)$. Then \Cref{thm: TClbub} give $2\leq \TC(\R H_{1,1})\leq 2$ and $2\leq \TC(\C H_{1,1})\leq 3$.  
Note that $\R H_{1,1}$ is a closed smooth manifold of dimension $1$. Therefore, $\R H_{1,1}\cong S^1$. Similarly, $\C H_{1,1}$ is a closed smooth manifold of dimension $2$. In fact, $\C H_{1,1}\cong S^2$. Therefore, $\TC(\C H_{1,1})=3$.

\item Suppose $(r,s)=(2,1)$. Then \Cref{thm: TClbub} give $3\leq \TC(\R H_{2,1})\leq 5$ and $3\leq \TC(\C H_{2,1})\leq 5$. Note that $\R H_{2,1}$ is a closed surface. Now observe that the fibre bundle in \eqref{eq: fbRHrs} looks like $\mathbb{R}P^{1} \stackrel{i}{\hookrightarrow} \mathbb{R}H_{2,1} \stackrel{p}{\longrightarrow} \mathbb{R}P^{1}.$ Therefore, $\R H_{2,1}$ can be either Torus or the Klein bottle. But from \Cref{thm: coring}, it follows that $\R H_{2,1}$ is the Klein bottle. Therefore, we get $\TC(\R H_{2,1})=5$ from \cite{TCKleinbottle}. 
One can see that $\C H_{2,1}\cong \C P^2\# \C P^2$. We have $5=\TC(\C P^2)\leq\TC(\C H_{2,1})$ using \cite[Theorem 2]{dranisSadyk}. Thus $\TC(\C H_{2,1})=5$.
\end{enumerate}

\end{example}

The following theorem relates $\TC(\R H_{r,s})$ and $\TC(\R P^s)+\TC(\R P^r)-1$.
\begin{theorem}
Let $s=2^{t_1}+1$ and $r=2^{t_2}$ with $t_2>1$ such that $s\leq r$. Then 
\begin{equation}\label{eq: prodineq}
 2(s+r-1)-1\leq \TC(\R H_{r,s}), \TC(\R P^r)+\TC( \R P^s)-1\leq 2(s+r-1).  
\end{equation} 
\end{theorem}
\begin{proof}
It follows from \eqref{eq: equaTCMil1} that,   $\TC(\R H_{r,s})\in \{2(s+r-1)-1,2(s+r-1)\}$, for $s=2^{t_1+1}$ and $r=2^{t_2}$ with $t_2>1$.  
On the contrary assume that 
\begin{equation}
  \TC(\R P^r)+\TC( \R P^s)-1< 2(s+r-1)-1.
\end{equation}
After putting $s=2^{t_1+1}$ and $r=2^{t_2}$ we get
\begin{equation}\label{eq: contra1}
\TC(\R P^{2^{t_2}})+\TC( \R P^{2^{t_1}+1})-1<2(2^{t_1}+2^{t_2})-1 .  
\end{equation}
From \cite[Corollary 8.2]{FarberTCproj} and \cite[Corollary 8.3]{FarberTCproj} we have $\TC(\R P^r)=2^{t_2+1}$ and $2^{t_1+1}\leq \TC(\R P^{2^{t_1}+1})$.
Therefore, 
\eqref{eq: contra1} gives 
\[2^{t_1+1}\leq \TC(\R P^s)< 2^{t_1+1}.\] This gives a contradiction. Thus the left inequality of \eqref{eq: prodineq} follows.

Note that the number of ones in diadic expansion of $s=2^{t_1}+1$ is $2$. Therefore,
it follows from \cite[Corollary 8.5]{FarberTCproj} that, $\TC(\R P^s)<2s $. 
Consequently, using \cite[Corollary 8.3]{FarberTCproj} we get \[\TC(\R P^{r})+\TC(\R P^s)-1< 2r+(2s-1)-1\leq2(s+r-1).\] This proves the right inequality of \eqref{eq: prodineq}.
\end{proof}

\begin{remark}
Recall the function $k(n)$ defined on the page $15$ of \cite{FarberTCproj}.
We observe that, if $k(s=2^{t_1}+1)>0$, then we have $\TC(\R P^s)<2s-1$. Consequently, we get \[\TC(\R P^r)+\TC(\R P^s)-1=2(s+r-1)-1\]when $s=2^{t_1}+1$ and $r=2^{t_2}$.  
\end{remark}

\begin{remark}
    T. Ganea in \cite{Ganeaconj} conjectured that for any finite CW-complex $X$ and $n>0$, $\ct(X\times S^n)=\ct(X)+\ct(S^n)-1$. The conjecture holds for a space $X$ whose cup-length coincides with the $\ct(X)-1$. Recall \Cref{prop: cat} shows $\cl(F H_{r,s})=\ct(F H_{r,s})-1$, where $F$ is $\R$ or $\C$. Therefore, Milnor manifolds satisfy the Ganea conjecture. 
\end{remark}

\section{Equivariant LS-category and equivariant topological complexity of Milnor manifolds}\label{sec: eq cat and tc}
Dey and Singh \cite{DeyMilnormfds} studied $\Z_2$ and $S^1$-actions on Milnor manifolds. They characterized Milnor manifolds which admit these actions under some conditions. 
In this section, we compute bounds on the equivariant LS-category and equivariant topological complexity of Milnor manifolds. 

\begin{proposition}\label{prop: eqcatRH}
Let $1< s< r$ and $r\not\equiv 2({\rm mod}~4)$. Suppose $\sigma$ is a fixed point free involution on $\R H_{r,s}$. Then
 \[(s+r)/2+1\leq \ct_{\left<\sigma\right>}(\R H_{r,s})\leq s+r.\]   
\end{proposition}
\begin{proof}
Recall that $\R H_{r,s}$ is a closed, smooth manifold of dimension $s+r-1$. Since $\sigma$ is a fixed point free involution on $\R H_{r,s}$, the space $\R H_{r,s}/\sigma$ is also a $(s+r-1)$-dimensional closed, smooth manifold. Since $\sigma$ is a fixed point free involution on $\R H_{r,s}$, we have $\ct_{\left<\sigma\right>}(\R H_{r,s})=\ct(\R H_{r,s}/\sigma)$ using \cite[Theorem 1.15]{Eqlscategory}. Therefore, using
\Cref{thm: catdimub}
we get the upper bound
\[\ct_{\left<\sigma\right>}(\R H_{r,s})\leq s+r.\]
Since $1<s<r$ with $r\not\equiv 2({\rm mod}~4)$, \Cref{thm: freeZ2 action} shows $\Z_2=\left<\sigma\right>$ acts trivially on the cohomology groups. Therefore, we can use the equivariant cohomology ring description of $\R H_{r,s}$ to compute the lower bound on $\ct_{\left<\sigma\right>}(\R H_{r,s})$.
Note that \Cref{thm: free invo milnor mfds} forces both $s$ and $r$ to be odd.
It follows from \Cref{thm: eq coho ring real milnor mfds}
that \[z\cdot x^{(s-1)/2}\cdot y^{(r-1)/2}\neq 0, \] where $|x|=2=|y|$ and $|z|=1$. In fact, $z\cdot x^{(s-1)/2}\cdot y^{(r-1)/2}$ generates $H^{s+r-1}(\R H_{r,s}/\sigma;\Z_2)$. Therefore,  $(s+r)/2\leq\cl(\R H_{r,s})$. We get $(s+r)/2+1\leq \ct_{\left<\sigma\right>}(\R H_{r,s})$ using \eqref{eq: cup-length}. This proves the proposition.
\end{proof}

\begin{remark}
 Observe that $\ct_{\left<\sigma\right>}(\R H_{1,1})=2$. On the other hand, we already know that $\R H_{1,1}\cong S^1$. Since $\R H_{1,1}/\left<\sigma\right>$ is a  smooth, closed manifold of dimension $1$, we have $\R H_{1,1}/\left<\sigma\right>\cong S^1$. Thus $\sigma$ has to be the antipodal involution.
\end{remark}

\begin{proposition}\label{prop: eq catCH}
Let $1< s< r$ and $r\not\equiv 2({\rm mod}~4)$. Let $\sigma$ be any fixed point free involution on $\C H_{r,s}$. Then
 \[(s+r)/2+2\leq \ct_{\left<\sigma\right>}(\C H_{r,s})\leq 2(s+r)-1.\]     
\end{proposition}
\begin{proof}
Recall that \Cref{thm: free invo milnor mfds} forces both $s$ and $r$ to be odd.
Then it follows from the \Cref{thm: eq coho ring complex milnor mfds} that \[z^2\cdot x^{(s-1)/2}\cdot y^{(r-1)/2}\neq 0,\] where $|z|=1$, $|x|=4=|y|$. Note that $|z^2\cdot x^{(s-1)/2}\cdot y^{(r-1)/2}|=2(s+r-1)$ and it generates $H^{2(s+r-1)}(\C H_{r,s})$.
Therefore, we get the desired lower bound using \eqref{eq: cup-length}. Note that $\C H_{r,s}/\sigma$ is a closed, smooth manifold of dimension $2(s+r-1)$. Then we get the upper bound from \Cref{thm: catdimub}. 
\end{proof}

We now describe the free $S^1$ action on $\R H_{r,s}$ when both $r$ and $s$ are odd . Note that only odd dimensional projective spaces admit free $S^1$-action. Let $r=2n+1$. Note that an element $w\in\R P^{r}$ can be described as a $[w_0,\dots,w_n]$, where $w_i\in 
\C$ for $0\leq i\leq n$. The $S^1$-action on $\R P^{r}$ is defined as follows: \[z\cdot [w_0,\dots,w_n]= [\sqrt{z}w_0,\dots,\sqrt{z}w_n] ,\] where $z\in S^1$. One can see that this $S^1$-action is free. Similarly, if $s=2m+1$, then we have a free $S^1$-action on $\R P^{s}$. Therefore, the diagonal action of $S^1$ on the product $\R P^r \times \R P^s$ is free. It can be seen that the real Milnor manifold $\R H_{r,s}$ is invariant under the diagonal $S^1$ action on $\R P^r\times \R P^s$. Thus, we have a free $S^1$ action on $\R H_{r,s}$. It was shown in \cite[Proposition 5.9]{DeyMilnormfds} that $\R H_{r,s}$ admit a free $S^1$ action if and only if both $r$ and $s$ are odd.
\begin{proposition}
Let $1\leq s\leq r$ such that both $r$ and $s$ are odd. Then 
\begin{equation}\label{eq: eq S1 cat }
(s+r)/2\leq \ct_{S^1}(\R H_{r,s})\leq s+r-1.    
\end{equation}
\end{proposition}
\begin{proof}
The hypothesis of the proposition shows that we can use equivariant cohomology ring information to compute the lower bound on $\ct_{S^1}(\R H_{r,s})$.
Observe that \Cref{thm: S1 eq coho ring} gives
$x^{(s-1)/2}\cdot y^{(r-1)/2}\neq 0$. In fact $x^{(s-1)/2}\cdot y^{(r-1)/2}$ generates $H^{s+r-2}(\R H_{r,s}/S^1;\Z_2)$. Then the left inequality of \eqref{eq: eq S1 cat } follows from \eqref{eq: cup-length}. Since both $r$ and $s$ are odd, it follows from \Cref{thm: free S1 action}  
that $\R H_{r,s}$ admit a free $S^1$-action. 
We have $\ct_{\left<\sigma\right>}(\R H_{r,s})=\ct(\R H_{r,s}/S^1)$ using \cite[Theorem 1.15]{Eqlscategory}. Note that $\R H_{r,s}/S^1$ is a closed, smooth manifold of dimension $s+r-2$. Therefore, the right inequality of \eqref{eq: eq S1 cat } follows from \Cref{thm: catdimub}.
\end{proof}

\begin{remark}
 Observe that if $(r,s)=(1,1)$, then $\ct_{S^1}(\R H_{1,1})=1$. Recall that $\R H_{1,1}\cong S^1$. So here we can say that the action of $S^1$ on itself is transitive.  
\end{remark}

We now prove the corresponding result for topological complexity.
\begin{proposition}\label{prop: eq tc FHrs}
Let $1< s< r$ and $r\not\equiv 2({\rm mod}~4)$ such that $s=2^{t_1}+1$, $r=2^{t_2}+1$. Let $\sigma$ be any fixed point free involution on $F H_{r,s}$. Then 
\[s+r-2\leq \TC_{\Z_2}(F H_{r,s}),\] where $F$ is either $\R$ or $\C$. 
\end{proposition}
\begin{proof}
Let $x\in H^{2}(\R H_{r,s}/\sigma; \Z_2)$ and $\bar{x}=1\otimes x+x\otimes 1$ be a zero-divisor. Then
\begin{align*}
\bar{x}^{s-2}&= \sum_{i=0}^{s-2}\binom{s-2}{i}x^i\otimes x^{s-i-2}\\
&= \binom{s-2}{(s-1)/2}x^{(s-1)/2}\otimes x^{(s-3)/2}+\binom{s-2}{(s-3)/2}x^{(s-3)/2}\otimes x^{(s-1)/2}\\
& = \binom{2^{t_1}-1}{2^{t_1-1}}x^{(s-1)/2}\otimes x^{(s-3)/2}+\binom{2^{t_1}-1}{2^{t_1-1}-1}x^{(s-3)/2}\otimes x^{(s-1)/2}\\
&= x^{(s-1)/2}\otimes x^{(s-3)/2}+x^{(s-3)/2}\otimes x^{(s-1)/2} \neq 0
\end{align*}
Similarly we can prove that $\bar{y}^{r-2}\neq 0$. Let $z\in H^{1}(\R H_{r,s}/\sigma;\Z_2)$. Let $A=x^{(s-1)/2}\cdot y^{(r-1)/2}$, $B=x^{(s-3)/2}\cdot y^{(r-3)/2}$, $C=x^{(s-3)/2}\cdot y^{(r-1)/2}$ and $D=x^{(s-1)/2}\cdot y^{(r-3)/2}$. Note that all these cohomology classes are non-zero. Therefore, the product
\begin{align*}
    \bar{z}\cdot \bar{x}^{s-2}\cdot \bar{y}^{r-2}& = A\otimes B\cdot z +C\otimes D\cdot z +D\otimes C\cdot z +B\otimes A\cdot z \\
    + & A\cdot z\otimes B +C\cdot z\otimes D +D\cdot z\otimes C +B\cdot z\otimes A
\end{align*}
is nonzero. Then the proposition follows from \cite[Theorem 5.15]{colmangranteqtc}. Similar techniques can be used to prove the proposition when $F=\C$.
\end{proof}

\begin{proposition}
Let $1\leq s\leq r$  such that $s=2^{t_1}+1$, $r=2^{t_2}+1$. Then 
\begin{equation}\label{eq: eq S1 tc }
s+r-2\leq \TC_{S^1}(\R H_{r,s}).    
\end{equation}
\end{proposition}
\begin{proof}
One can use the $S^1$-equivariant cohomology ring description from \Cref{thm: S1 eq coho ring} and follow the similar steps as in \Cref{prop: eq tc FHrs} to prove $\bar{w}\cdot \bar{x}^{s-2}\cdot \bar{y}^{r-2}\neq 0$, when $s=2^{t_1}+1$, $r=2^{t_2}+1.$ Thus the proposition follows from \cite[Theorem 5.15]{colmangranteqtc}.
\end{proof}

Since $\R H_{r,s}$ is free metrizable space, \cite[Proposition 1.15]{Eqlscategory} gives $\ct_{\Z_2}(\R H_{r,s}\times \R H_{r,s})=\ct(\R H_{r,s}\times \R H_{r,s}/\Z_2)$.
Therefore, using \Cref{thm: catdimub} we get 
\[\TC_{\Z_2}(\R H_{r,s})\leq \ct((\R H_{r,s}\times \R H_{r,s})/\Z_2)\leq  2\cdot \
\mathrm{dim}(\R H_{r,s})+1.\]
Therefore, the following proposition follows from \Cref{thm: TClbub}.

\begin{proposition}
If $s=2^{t_1}+1$ and $r=2^{t_2}$, then \[2(s+r-1)-1\leq \TC(\R H_{r,s})\leq \TC_{\Z_2}(\R H_{r,s})\leq 2(s+r-1)+1.\]
\end{proposition}

\section{Generalized projective product spaces with one of the components as Milnor manifolds}\label{sec: cat and tc gpps}
Sarkar and Zvengrowski \cite{sarkargpps} defined the generalized projective product spaces and studied topological aspects of some of its classes. The author of this paper and Sarkar \cite{DaundSarkargpps}, obtained upper bounds on the LS-category and topological complexity of generalized projective product spaces (see \cite[Proposition 2.4, Proposition 2.5]{DaundSarkargpps}). 
In this section, we compute bounds on the LS-category and the topological complexity of generalized projective product spaces with Milnor manifolds as one of its components.

Under certain conditions, authors of \cite{DeyMilnormfds} classified Milnor manifolds that admit free involutions. 
    Let $1 < s < r$ and $r\not\equiv 2({\rm mod}~ 4)$. Then, it was shown in \cite{DeyMilnormfds} that both $\R H_{r,s}$ and $\C H_{r,s}$ admit free involution if and only if both $r$ and $s$ is
odd. We denote such free involution on $\R H_{r,s}$ by $\sigma_{\R}$ and on $\C H_{r,s}$ by $\sigma_{\C}$.
Let $N$ be a path connected topological space with a fixed point free involution $\tau$. Define 
\[X(\R H_{r,s},N):=\displaystyle \frac{\R H_{r,s}\times N}{(x,y)\sim (\sigma_{\R}(x),\tau(y))},\] and \[X(\C H_{r,s},N):=\displaystyle \frac{\C H_{r,s}\times N}{(x,y)\sim (\sigma_{\C}(x),\tau(y))},\]
Note that we have fibre bundles 
\begin{equation}\label{eq: fbRH}
\R H_{r,s}\xhookrightarrow{}X(\R H_{r,s},N)\stackrel{\mathfrak{p}}\longrightarrow N/\tau    
\end{equation}

and
\begin{equation}\label{eq: fbCH}
\C H_{r,s}\xhookrightarrow{}X( \C H_{r,s},N)\stackrel{\mathfrak{q}}\longrightarrow N/\tau.    
\end{equation}

\begin{proposition}\label{prop: catXRHN}
Let $0 \leq s \leq r$. Let $N$ be the path connected space with a fixed point free involution $\tau$. Then
\[\ct(X(\R H_{r,s}, N))\leq \ct(N/\tau) +s+r-1.\]\
\end{proposition}
\begin{proof}
Recall the fibre bundle in \eqref{eq: fbRH}. Then using \cite[Proposition 2.4]{DaundSarkargpps}, we get $\ct(X(\R H_{r,s}, N))\leq \ct(N/\tau) + \ct_{\left<\sigma\right>}(R H_{r,s})-1$. Then the proposition follows from \Cref{prop: eqcatRH}.
\end{proof}

The following result follows from the \cite[Proposition 2.5]{DaundSarkargpps}.
\begin{proposition}
Let $0 \leq s \leq r$. Let $N$ be the path connected space with a fixed point free involution $\tau$ and $\sigma$ be a free nvolution on $\R H_{r,s}$. 
Then \[\TC(X(\R H_{r,s}, N))\leq \ct(N/\tau\times N/\tau)+ \TC_{\left<\sigma\right>}(\R H_{r,s})-1\]   
\end{proposition}

We now compute the mod-$2$ cohomology ring of $X(\C H_{r,s}, N)$, where $N$ is a simply connected space.

\begin{proposition}\label{prop: coho ring XCH}
Let $1 < s < r$ be odd integers such that $r\not\equiv 2({\rm mod}~ 4)$ and $N$ be a simply connected, path connected topological space with a fixed point free involution $\tau$. Then \[H^{\ast}(X(\C H_{r,s}, N);\Z_2)=H^{\ast}(\C H_{r,s};\Z_2)\otimes H^{\ast}(N/\tau;\Z_2).\]   
\end{proposition}
\begin{proof}
Note that $\C H_{r,s}$ is a compact, simply connected, and path connected space. Since $1 < s < r$ be odd integers such that $r\not\equiv 2({\rm mod}~ 4)$, using \Cref{thm: free invo milnor mfds} we have a free involution $\sigma_{\C}$ on $\C H_{r,s}$. Then it follows from \Cref{thm: freeZ2 action}
that 
$\sigma_{\C}^{\ast}$ is trivial. Therefore, the proposition follows from \cite[Theorem 4.11]{DaundSarkargpps}.
\end{proof}

We use the mod-$2$ cohomology ring information to compute bounds on the LS-category and topological complexity of $X(\C H_{r,s})$.

\begin{proposition}\label{prop: cat XCHrsN}
Let $1 < s < r$ be odd integers such that $r\not\equiv 2({\rm mod}~ 4)$ and $N$ be the simply connected path connected space. Then
\begin{equation}\label{eq: catXCH}
r+s+\cl_{\Z_2}(N/\tau)\leq \ct(X(\C H_{r,s}, N))\leq \ct(N/\tau) +2(s+r))-2. 
\end{equation}
\end{proposition}
\begin{proof}
It follows from the \Cref{prop: coho ring XCH} that \[\cl_{\Z_2}(X(\C H_{r,s},N))=r+s+\cl_{\Z_2}(N/\tau)-1.\] Therefore, the  left inequality of \eqref{eq: catXCH} follows from \eqref{eq: cup-length}. Using \cite[Proposition 2.4]{DaundSarkargpps} we get 
\[\ct(X(\C H_{r,s}, N))\leq \ct(N/\tau) +\ct_{\left<\sigma_{\C}\right>}(\C H_{r,s})-1.\]
Now the right inequality of \eqref{eq: catXCH} follows from \Cref{prop: eq catCH}.
\end{proof}

We now prove the corresponding result for the topological complexity.
Let $\sigma_{\C}$ be a fixed point free involution on $\C H_{r,s}$ described at the beginning of this section.
\begin{proposition}
Let $1 < s < r$ be an odd integers such that $r\not\equiv 2({\rm mod}~ 4)$ and $N$ be a simply connected path connected space. Then 
\begin{equation}\label{prop: TCXCH}
\zl_{\Z_2}(N/\tau)+\zl_{\Z_2}(\C H_{r,s})+1\leq \TC(X(\C H_{r,s}, N))\leq \ct(N/\tau\times N/\tau)+ \TC_{\left<\sigma_{\C}\right>}(\C H_{r,s})-1.    
\end{equation}   
\end{proposition}
\begin{proof}
The upper bound on $\TC(X(\C H_{r,s}, N)$ follows from the \cite[Proposition 2.5]{DaundSarkargpps} and the lower bound follows from the \Cref{prop: coho ring XCH} and zero-divisor-cup-length calculations.
\end{proof}

The unit $n$-sphere in $\R^{n+1}$ is defined as  \[S^n :=\{(u_1, \ldots, u_{n+1})\in \R^{n+1} ~|~ u_1^2 + \cdots + u_{n+1}^2 = 1\}.\] 
In \cite{Davis}, Davis discovered new topological spaces called projective product spaces and studied their topological properties. Let $(n_1,\dots,n_k)$ be a tuple of non-negative integers. Define
\begin{equation}\label{eq_ppsp}
P(n_1, \ldots, n_k):= \frac{S^{n_1}\times \cdots \times S^{n_k}}{({\bf x}_1, \dots, {\bf x}_{k})\sim (-{\bf x}_1, \dots, -{\bf x}_{k})}.
\end{equation}
The identification space $P(n_1, \ldots, n_k)$ is called a \emph{projective product spcace} corresponding to the tuple $(n_1, \ldots, n_k)$. Observe that if $k=1$, then $P(n_1)=\R P^{n_1}$.

Consider the following identification space: 
\[X(FH_{r,s},n_1,\dots,n_k)=\frac{FH_{r,s}\times S^{n_1}\times \cdots \times S^{n_k}}{(z,{\bf x}_1, \dots, {\bf x}_{k})\sim (\sigma _F(z),-{\bf x}_1, \dots, -{\bf x}_{k})},\] where $F$ is $\R$ or $\C$ and $\sigma_F$ is an involution on $FH_{r,s}$. 
Observe that, we have a fibre bundle 
\begin{equation}\label{eq: fb FHXFHP}
FH_{r,s}\xhookrightarrow{}X(FH_{r,s}, n_1,\dots,n_k)\xrightarrow{\mathfrak{p}}P(n_1,\dots,n_k),    
\end{equation}
where $P(n_1,\dots,n_k)$ is the projective product space defined in \eqref{eq_ppsp}. 

The following result is a straightforward consequence of \Cref{prop: catXRHN} and \cite[Theorem 1.2]{SeherVandem}.
\begin{proposition}
Let $1 < s < r$ be an odd integers such that $r\not\equiv 2({\rm mod}~ 4)$ and $2\leq n_1 \leq \cdots \leq n_k$.  
Then \[\ct(X(\R H_{r,s}, n_1,\dots,n_k))\leq n_1+k+s+r-1.\]    
\end{proposition}

We give a description of the mod-$2$ cohomology ring of $X(\C H_{r,s},n_1,\dots,n_k)$.
\begin{proposition}\label{prop: cohoXCHP}
 Let $1 < s < r$ be an odd integers such that $r\not\equiv 2({\rm mod}~ 4)$ and $2\leq n_1 \leq \cdots \leq n_k$. Then  \[H^*(X(\C H_{r,s}, n_1, \ldots, n_k)); \Z_2) \cong H^*(\R P^{n_1};\Z_2)\otimes \Lambda(\beta_2,\ldots, \beta_k) \otimes H^*(\C H_{r,s}; \Z_2),\] where $\Lambda(-)$ is an exterior algebra over $\Z_2$.
\end{proposition}
\begin{proof}
The proof follows from   \Cref{prop: coho ring XCH} and \cite[Theorem 2.1]{Davis}.   
\end{proof}

\begin{proposition}
Let $1 < s < r$ be an odd integers such that $r\not\equiv 2({\rm mod}~ 4)$ and $2\leq n_1\leq \dots \leq n_r$. Then
\begin{equation}\label{eq: catXCHP}
r+s+n_1+k-1\leq \ct(X(\C H_{r,s}, n_1,\dots, n_k))\leq n_1+k +2(s+r)-2. 
\end{equation}
\end{proposition}
\begin{proof}
It follows from the \Cref{prop: cohoXCHP} and \cite[Proposition 2.3]{SeherVandem} that \[\cl_{\Z_2}(X(\C H_{r,s},n_1,\dots,n_k))=r+s+n_1+k-2.\] Therefore, the left inequality of \eqref{eq: catXCHP} follows from \eqref{eq: cup-length}.  
Using \cite[Proposition 2.4]{DaundSarkargpps}, we get
\[\ct(X(\C H_{r,s}, n_1,\dots, n_k))\leq n_1+k +\ct_{\left<\sigma_{\C}\right>}(\C H_{r,s})-1.\]
Then the right inequality of \eqref{eq: catXCHP} follows from and \cite[Theorem 1.2]{SeherVandem} and \Cref{prop: eq catCH}.
\end{proof}

We now prove the corresponding result for the topological complexity.

\begin{proposition}
Let $1 < s < r$ be an odd integers such that $r\not\equiv 2({\rm mod}~ 4)$ and $2\leq n_1\leq \dots \leq n_r$. Then 
\begin{equation}\label{eq: TCXCHP}
\zl_{\Z_2}(\R P^{n_1})+k+\zl_{\Z_2}(\C H_{r,s})\leq \TC(X(\C H_{r,s},n_1,\dots,n_k))\leq 2(n_1+k-1)+ \TC_{\left<\sigma_{\C}\right>}(\C H_{r,s}).    
\end{equation}   
\end{proposition}
\begin{proof}
One can see that the $\zl(P(n_1,\dots,n_k))=\zl(\R P^{n_1})+k-1$ (see \cite[Page 2]{SeherVandem}). 
Then the left inequality of \eqref{eq: TCXCHP} follows from the zero-divisor-cup-length calculations using \Cref{prop: cohoXCHP}. Now recall the fibre bundle \eqref{eq: fb FHXFHP}. Then we get the inequality 
\[\TC(X(\C H_{r,s},n_1,\dots,n_k))\leq cat(P(n_1,\dots,n_k)\times P(n_1,\dots,n_k))+\TC_{\left<\sigma_{\C}\right>}(\C H_{r,s})-1\] 
using \cite[Proposition 2.5]{DaundSarkargpps}. Now the right inequality of \eqref{eq: TCXCHP} follows from \cite[Theorem 1.2]{SeherVandem} and \cite[Proposition 1.37]{CLOT}.

\end{proof}

\begin{discussion}\label{discussion}
\normalfont{
Observe that if $s=r$, then the mapping \[\big([z_0,\dots,z_s],[w_0,\dots,w_s]\big) \to \big([w_0,\dots,w_s],[z_0,\dots,z_s]\big)\] define a fixed point free involution on $\C H_{s,s}$. We denote this involution by $I_{\C}$. We can have a  similar fixed point free involution on $\R H_{s,s}$. Let us denote it by $I_{\R}$. 

We now construct a $I_{\C}$-invariant categorical cover of $\C H_{r,s}$. Let $U_i=\{[z_0,\dots,z_s]\in \C^s ~\mid~ z_i\neq 0\}$ for $0\leq i\leq s$. Then note that $U_i\cong \C^s$ and $\{U_i\mid 0\leq i\leq s\}$ forms a categorical cover of $\C P^s$. Consider the open sets of $\C P^s\times \C P^s$ given by $W_j=\cup_{t=0}^{j}U_t\times U_{j-t}$
for $1\leq j \leq 2s+1$. Note that each $W_j$ is $I_{\C}$-invariant and the collection 
\begin{equation}\label{eq: cat cover of CH}
 \{W_j \mid 1\leq j\leq 2s+1\}   
\end{equation}
forms a categorical cover of $\C P^s\times \C P^s$. Observe that for each $0\leq t\leq j$, \[\C H_{s,s}\cap U_t\times U_{j-t}=\bigg\{\big([z_0,\dots,z_s],[w_0,\dots,w_s]\big)\in (\C P^s)^2 \mid z_t\neq 0, z_{j-t}\neq 0 ~\&~ \sum_{i=0}^{s}z_i\overline{w}_i=0\bigg\}.\]
One can see that $\C H_{s,s}\cap U_t\times U_{j-t}$ is homeomorphic to a hyperplane in $\C^s\times \C^s$. Consequently, $\C H_{s,s}\cap U_t\times U_{j-t}$ is contractible for $1\leq t\leq j$. This implies the sets $W_j$'s are categorical in $\C H_{s,s}$ for each $1\leq j\leq 2s+1$. Thus, the collection $\{\C H_{s,s}\cap W_j \mid 1\leq j\leq 2s+1\}$ forms an $I_{\C}$-invariant categorical cover of $\C H_{s,s}$. Similar computations hold for $\R H_{s,s}$.}
\end{discussion}

\begin{proposition}
 Let $s> 0$ and $N$ be a simply connected, path connected space with a fixed point free involution $\tau$ . Then 
\begin{equation}\label{eq: catssXCH}
    \ct(X(F H_{s,s},N))\leq \ct(N/\tau) + 2s,
\end{equation}
where $F$ is either $\R$ or $\C$.
\end{proposition}
\begin{proof}
Using \cite[Proposition 2.4]{DaundSarkargpps} to we have  \[\ct(X(\C H_{s,s},N))\leq \ct(N/\tau)+q-1,\] where $q$ is the number of $I_{\C}$-invarinat categorical open sets which covers $\C H_{r,s}$. From  the \Cref{discussion}, we have $q\leq 2s+1$. Now the inequality in \eqref{eq: catssXCH} follows.
 Similarly, we get $\ct(X(\R H_{s,s},N))\leq \ct(N/\tau) + 2s$.
 This proves the proposition.
\end{proof}

Now we do similar computations for topological complexity.

\begin{proposition}
Let $s>0$  and $N$ be a simply connected, path connected space with a fixed point free involution $\tau$ . Then 
\begin{equation}\label{eq: tcssXCH}
    \TC(X(F H_{s,s},N))\leq \ct(N/\tau\times N/\tau) + 4s,
\end{equation}
where $F$ is either $\R$ or $\C$.
\end{proposition}
\begin{proof}
Using \cite[Proposition 2.5]{DaundSarkargpps} we get \[\TC(X(\C H_{s,s},N))\leq \ct(N/\tau\times N/\tau) +q-1,\] where $q$ is the number of $(I_{\C}\times I_{\C})$-invariant open sets on which the free path space fibration $\pi: P(\C H_{s,s})\to \C H_{s,s}\times \C H_{r,s}$ has continuous sections. Now we can use \eqref{eq: cat cover of CH} and \eqref{eq_ls_cat_prod} to construct   $(I_{\C}\times I_{\C})$-invariant categorical cover of $\C H_{s,s}\times \C H_{s,s}$ with $4s+1$ many open sets. Since the such open cover of $\C H_{s,s}\times \C H_{s,s}$ is categorical, we have local continuous sections of $\pi$. This proves the right inequality of \eqref{eq: tcssXCH}.
Similarly, we can prove that $\TC(X(\R H_{s,s},N))\leq \ct(N/\tau\times N/\tau) + 4s$.
\end{proof}

\begin{remark}
Observe that if $2\leq n_1\leq \dots\leq n_k$ and $N=S^{n_1}\times \dots\times S^{n_k}$ such that $N/\tau=P(n_1,\dots,n_k)$, then using \cite[Theorem 1.2, Proposition 2.3]{SeherVandem} and \eqref{eq: catssXCH} we get \[\ct(X(F H_{s,s},N))\leq n_1+k + 2s,\]   and using \eqref{eq: tcssXCH} we get \[\TC(X(F H_{s,s},N))\leq 2(n_1+k) + 4s-1,\] where $F$ is either $\R$ or $\C$. 
\end{remark}

Next, we define another class of generalized projective product spaces related to Milnor manifolds and compute bounds on LS-category and topological complexity of these spaces.
Let $N=S^{n_1}\times \cdots \times S^{n_k}$ be the product of unit spheres with the involution $\tau$ determined by the antipodal involution on each factor $S^{n_j}$ for $1\leq j\leq k$. 

Now we recall the involution $\tau_j$ on $S^{n_j}$ defined in \cite[Example 3.1]{sarkargpps}.
\begin{equation}\label{eq: invo prodsphere}
  \tau_j( (y_1, \ldots, y_{p_j}, y_{p_j+1}, \ldots, y_{n_j+1}) ):= (y_1, \ldots, y_{p_j}, -y_{p_j+1}, \ldots, -y_{n_j+1}),  
\end{equation}
for some $0 \leq p_j \leq n_j$ and $1\leq j\leq k$.
 Then we have $\Z_2$-action on the product $S^{n_1}\times \dots \times S^{n_{k}}$ via product of these involutions
 \begin{equation}\label{eq: invo on prod sphere}
  \tau_1\times \dots \times \tau_{k}: S^{n_1}\times \dots \times S^{n_{k}}\to S^{n_1}\times \dots \times S^{n_{k}}   
 \end{equation}
We denote this involution by $\mathfrak{I}$.  It can be easily seen that if all $p_j=0$ then $\tau_j$ is the antipodal involution on $S^{n_j}$ and if $p_j=n_j$, then $\tau_j$ is the reflection across the hyperplane $y_{n_j+1}=0$ in the Euclidean space $\R^{n_j+1}$.

Let $s$ and $r$ be an odd integers such that $1 < s < r$ and $r\not\equiv 2({\rm mod}~ 4)$. Then by \cite[Corollary 5.6 and Corollary 5.8 ]{DeyMilnormfds} both $\R H_{r,s}$ and $\C H_{r,s}$ admit free involutions. 
Consider $\R H_{r,s}$ with a free involution $\sigma_{\R}$ and $\C H_{r,s}$ with $\sigma_{\C}$. 
We now consider the identification space defined as follows: 
\begin{equation}\label{eq:prodsphere_RH}
X((n_1, p_1), \ldots, (n_k,p_k), \R H_{r,s}) :=\frac{ S^{n_1}\times \cdots \times S^{n_k}\times \R H_{r,s}}{(x_1,\dots ,x_{k}, y)\sim (\tau_1(x_1),\dots ,\tau_k(x_k), \sigma_{\R}(y))},
\end{equation} and 
\begin{equation}\label{eq:prodsphere_CH}
X((n_1, p_1), \ldots, (n_k,p_k), \C H_{r,s}) :=\frac{ S^{n_1}\times \cdots \times S^{n_k}\times \C H_{r,s}}{(x_1,\dots ,x_{k}, y)\sim (\tau_1(x_1),\dots ,\tau_k(x_k), \sigma_{\C}(y))},
\end{equation}
where $\tau_j$ is a reflection defined as in \eqref{eq: invo prodsphere} for $1\leq j\leq k$. Note that the spaces defined in \eqref{eq:prodsphere_CH} and \eqref{eq:prodsphere_RH} are  generalized projective product spaces. 
Consequently, we have a fibre bundle 
\begin{equation}\label{eq: fibundXppsFH}
    S^{n_1}\times\dots\times S^{n_k}\xhookrightarrow{} X((n_1, p_1), \ldots, (n_k,p_k), F H_{r,s})\stackrel{\mathfrak{p}}\longrightarrow F H_{r,s}/\sigma_F,.
\end{equation}
where $F$ is either $\R$ or $\C$ and $\mathfrak{p}$ is the map induced by the orbit map $F H_{r,s}\to F H_{r,s}/\sigma_F$.

Let $F$ be either $\R$ or $\C$ and $\alpha_F$ be the first Stiefel-Whitney class of the line bundle associated with the double cover $F H_{r,s} \to F H_{r,s}/\sigma_{F}$. We denote the mod-$2$ exterior algebra by $\Lambda(-)$ and the total Steenrod square by ${\rm Sq} = \sum_{n \geq 0} {\rm Sq}^n$.
Now we are ready to describe the mod-$2$ cohomology ring spaces defined in \eqref{eq:prodsphere_RH} and \eqref{eq:prodsphere_CH}.
\begin{proposition}\label{prop: coho ring XPPSRH}
Let $1 < s < r$ be an odd integers such that $r\not\equiv 2({\rm mod}~ 4)$. Let $n_1 \leq \cdots \leq n_k$ and $F$ is either $\R$ or $\C$ . Then  
\[H^*(X((n_1, p_1), \ldots, (n_k,p_k), F H_{r,s});\Z_2)\cong H^*(F H_{r,s}/\sigma_{F}, \Z_2) \otimes \Lambda(\beta_1,\ldots, \beta_k),\]

where  $|\beta_j|=n_j,$ ${\rm Sq}(\beta_j)=(1+\alpha)^{n_j+1-p_j} \beta_j$, and and $p_j \geq 1$ for $1 \leq j \leq k$.   
\end{proposition}
\begin{proof}
Follows from \cite[Theorem 4.]{DaundSarkargpps}. 
\end{proof}

\begin{proposition}\label{prop: cat XPPSRH}
Let $1 \leq p_j \leq n_j$ and $1 < s < r$ be an odd integers such that $r\not\equiv 2({\rm mod}~ 4)$. Then,
we have the following inequalities. 
\begin{equation}\label{eq: cat XPPSRH }
 (s+r)/2+k+1\leq \mathrm{cat}(X((n_1,p_1) \ldots, (n_k,p_k), \R H_{r,s}))\leq r+s+k\end{equation}
\end{proposition}
\begin{proof}
Let $\sigma_{\R}$ be a fixed point free involution on $\R H_{r,s}$.
 Using \Cref{prop: coho ring XPPSRH}, we have \[\cl_{\Z_2}(X((n_1,p_1) \ldots, (n_k,p_k), \R H_{r,s}))=\cl_{\Z_2}(\R H_{r,s}/\left<\sigma_{\R}\right>)+ k.\] 
 It can be seen that the cup-length of $\R H_{r,s}/\left<\sigma\right>$ is $(s+r)/2$ (see proof of \Cref{prop: eqcatRH}).
 Then we get the left inequality of \eqref{eq: cat XPPSRH } using \eqref{eq: cup-length}.
Note that the \cite[Proposition 2.4]{DaundSarkargpps} give
\[\mathrm{cat}(X((n_1,p_1) \ldots, (n_k,p_k), \R H_{r,s}))\leq \mathrm{cat}(\R H_{r,s}/\sigma_{\R})+q-1,\] where $q$ is the number of $\mathfrak{I}$-invariant open categorical open sets that covers $\prod_{i=1}^{k}S^{n_i}$. Consider $U_{i_1}=S^{n_i}\setminus \{e_{n_i}(1)\}$ and $U_{i_2}=S^{n_i}\setminus \{-e_{n_i}(1)\}$, where $e_{n_i}(1)=(1,0,\dots,0)\in S^{n_i}$. Observe that $\tau_i(U_{i_1})=U_{i_1}$ and $\tau_i(U_{i_2})=U_{i_2}$. Therefore, $\{U_{i_1}, U_{i_2}\}$ is an $\left<\tau_i\right>$-invariant open cover of $S^{n_i}$. One can see that they are also $\left<\tau_i \right>$-categorical for $i=1, \ldots, k$. 
Now using \cite[Proposition 1.37]{CLOT}, we have $\ct(\prod_{i=1}^{k}S^{n_i}) \leq k+1$. Thus, one can construct a categorical open cover of $\prod_{i=1}^{k}S^{n_i}$ with $k+1$ open sets invariant under the $(\tau_1\times \dots\times \tau_k)$-action. Now right inequality of  \eqref{eq: cat XPPSRH } follows from \Cref{prop: eqcatRH}. 
\end{proof}

\begin{remark}
Note that the open cover constructed in \Cref{prop: cat XPPSRH} appeared in the proof of \cite[Proposition 4.3]{DaundSarkargpps}.  
\end{remark}

\begin{proposition}
Let $1 \leq p_j \leq n_j$ for $1\leq j\leq k$ and  $1 < s < r$ be an odd integers such that $r\not\equiv 2({\rm mod}~ 4)$.   Then,
we have the following inequalities. 
\begin{equation}\label{eq: catXPPSCHrs}
(s+r)/2+k+2\leq \mathrm{cat}(X((n_1,p_1) \ldots, (n_k,p_k), \C H_{r,s}))\leq 2(s+r)+k-1.    
\end{equation}
 
\end{proposition}
\begin{proof}
Let $\sigma_{\C}$ be a fixed point free involution on $\C H_{r,s}$.
From \Cref{prop: coho ring XPPSRH}, we get \[\cl_{\Z_2}(X((n_1,p_1) \ldots, (n_k,p_k), \C H_{r,s}))=\cl_{\Z_2}(\C H_{r,s}/\sigma_{\C})+k.\] One can see that  $(s+r)/2+1\leq \cl_{\Z_2}(\C H_{r,s}/\sigma_{\C})$ from \Cref{thm: eq coho ring complex milnor mfds}. Therefore, we get the left inequality of \eqref{eq: catXPPSCHrs}. Consider the fibre bundle in\eqref{eq: fibundXppsFH}. Then we get the following inequality follows from \cite[Proposition 2.4]{DaundSarkargpps}
    \[\mathrm{cat}(X((n_1,p_1) \ldots, (n_k,p_k), \C H_{r,s}))\leq \mathrm{cat}(\C H_{r,s}/\sigma_{\C})+k.\] 
Now the right inequality of \eqref{eq: catXPPSCHrs} follows from \Cref{prop: eq catCH}.    
\end{proof}

\begin{proposition}\label{prop: tcxppsrh}
Let $1 < s < r$ be an odd integers such that $r\not\equiv 2({\rm mod}~ 4)$ and $2\leq n_1 \leq \cdots \leq n_k$. Suppose $p_j >1$ for $j=1, \ldots, k$.
Then  
\begin{equation}\label{eq: tcXPPSRH}
   \zl_{\Z_2}(\R H_{r,s}/\sigma_{\R})+ k+1\leq\TC(X((n_1,p_1),\dots,(n_k,p_k),\R H_{r,s})\leq 2(s+r+k)-1.  
\end{equation}
Moreover, if $r=2^{t_1}+1$ and $s=2^{t_2}+1$, then \begin{equation}\label{eq: tlbcXPPSRH}
 s+r+ k-1\leq\TC(X((n_1,p_1),\dots,(n_k,p_k),\R H_{r,s})\leq 2(s+r+k)-1.  
\end{equation}
\end{proposition}
\begin{proof}
Let $\sigma_{\R}$ be a fixed point free involution on $\R H_{r,s}$.
Using \Cref{prop: coho ring XPPSRH}, we have \[\zl_{\Z_2}(X((n_1,p_1), \ldots, (n_k,p_k), \R H_{r,s}))=\zl_{\Z_2}(\R H_{r,s}/\left<\sigma_{\R}\right>)+ k.\] This give us the left inequality of \eqref{eq: tcXPPSRH}. Now consider the fibre bundle described in \eqref{eq: fibundXppsFH}. Then using \cite[Proposition 2.5]{DaundSarkargpps} we have
\[\TC(X((n_1,p_1)\dots,(n_k,p_k),\R H_{r,s})\leq \ct(\R H_{r,s}/\left<\sigma_{\R}\right>\times  \R H_{r,s}/\left<\sigma_{\R}\right>)+ q-1, \]
 where $q$ is the number of $(\mathfrak{I}\times \mathfrak{I})$-invariant open sets of $\prod_{i=1}^{k}S^{n_i}\times \prod_{i=1}^{k}S^{n_i}$
 which admit a local continuous sections of $\pi: P(\prod_{i=1}^{k}S^{n_i})\to \prod_{i=1}^{k}S^{n_i}\times \prod_{i=1}^{k}S^{n_i}$.  
 Recall that we constructed $\mathfrak{I}$-invariant categorical open cover of $\prod_{i=1}^{k}S^{n_i}$ with $k+1$ many open sets in \Cref{prop: cat XPPSRH}. Therefore, using \eqref{eq_ls_cat_prod} we can get the $(\mathfrak{I}\times \mathfrak{I})$-invariant categorical  open cover of  $\prod_{i=1}^{k}S^{n_i}\times \prod_{i=1}^{k}S^{n_i}$ with $2k+1$ many open sets. Now, this open cover serve our purpose. Now, we get the right inequality of \eqref{eq: tcXPPSRH} using \Cref{prop: eqcatRH}.

 The left inequality in \eqref{eq: tlbcXPPSRH} follows from the calculation of $\zl_{\Z_2}(\R H_{r,s}/\sigma_{\R})$ as done in \Cref{prop: eq tc FHrs}.
\end{proof}

We now prove a similar result for the topological complexity.
\begin{proposition}
Let $1 < s < r$ be an odd integers such that $r\not\equiv 2({\rm mod}~ 4)$ and $2\leq n_1 \leq \cdots \leq n_k$. Suppose $p_j >1$ for $j=1, \ldots, k$.
Then  
\begin{equation}\label{eq: tcXPPSCH}
  \zl_{\Z_2}(\C H_{r,s}/\sigma_{\C })+ k+1\leq\TC(X((n_1,p_1)\dots,(n_k,p_k),\C H_{r,s})\leq 4(s+r+k)-2k-3. 
\end{equation}
Moreover, if $r=2^{t_1}+1$ and $s=2^{t_2}+1$, then 
\begin{equation}\label{eq: lbtcXPPSCH}
  s+r+ k-1\leq\TC(X((n_1,p_1)\dots,(n_k,p_k),\C H_{r,s})\leq 4(s+r+k)-2k-3. 
\end{equation}
\end{proposition}
\begin{proof}
The proof is similar to that of \Cref{prop: tcxppsrh}.
\end{proof}

\section{Acknowledgement}
The author would like to thank the reviewer for careful reading and for valuable suggestions to improve the exposition of the article. The author also would like to thank Goutam Mukherjee, Soumen Sarkar and Bittu Singh for helpful discussion's.

\bibliographystyle{plain} 
\bibliography{references}

\end{document}